\documentclass[11pt]{amsart}

\usepackage{amssymb}

\newcommand{\mylabel}[1]{\mbox{\rm \quad L(#1)}\label{#1}}

\def\proof{\par\medbreak\noindent{\bf Proof}\quad\enspace}
\def\note{\par\medbreak\noindent{\bf Note}\quad\enspace}
\def\remark{\par\medbreak\noindent{\bf Remark}\quad\enspace}
\def\example{\par\medbreak\noindent{\bf Example}\quad\enspace}

\def\platz{\par\medbreak}

   {\end{trivlist}\platz}

\def\qed{\vbox{\hrule
  \hbox{\vrule\hbox to 5pt{\vbox to 8pt{\vfil}\hfil}\vrule}\hrule}}

\newtheorem{@theorem}{Theorem}[section]
\newtheorem{@prop}[@theorem]{Proposition}
\newtheorem{@lemma}[@theorem]{Lemma}
\newtheorem{@cor}[@theorem]{Corollary}
\newtheorem{@conjecture}[@theorem]{Conjecture}
\newtheorem{@dfn}[@theorem]{Definition}

\newenvironment{theorem}[1]{\begin{@theorem}\mylabel{#1}}{\end{@theorem}}
\newenvironment{prop}[1]{\begin{@prop}\mylabel{#1}}{\end{@prop}}
\newenvironment{lemma}[1]{\begin{@lemma}\mylabel{#1}}{\end{@lemma}}

\newtheorem{@satz}[@theorem]{Satz}
\newtheorem{@satzdef}[@theorem]{Satz und Definition}
\newtheorem{@kor}[@theorem]{Korollar}
\newtheorem{@folg}[@theorem]{Folgerung}
\newtheorem{@problem}[@theorem]{Problem}
\newtheorem{@notiz}[@theorem]{Notiz}

\newcommand {\eps}{\varepsilon}

\newcommand\const{{\rm const}}

\newcommand\dist{{\rm dist\,}}

\renewcommand\Im{{\rm Im\,}}

\renewcommand\Re{{\rm Re\,}}

\newcommand\supp{{\rm supp\,}}

\newcommand {\R}{\mathbb{R}}

\newcommand {\C}{\mathbb{C}}
\newcommand {\N}{\mathbb{N}}
\newcommand {\nnull}{\N_0}
\newcommand {\Q}{\mathbb{Q}}

\renewcommand{\mylabel}[1]{\label{#1}}

\newcommand {\qi} {\approxeq}  
\newcommand {\gtilde} {{\tilde{g}}}

\newcommand {\ptilde} {{\tilde{p}}}
\newcommand {\qtilde} {{\tilde{q}}}
\newcommand {\rhotilde} {{\tilde{\rho}}}

\newcommand {\lt} {\tilde{l}}

\newcommand {\gSN} {{g_{S^{N-1}}}}
\newcommand {\pij} {\pi_{(j)}}
\newcommand {\Ubar} {\overline{U}}

\newcommand {\DGB}{D_{\text{GB}}}
\newcommand {\dom}{\mathcal{D}}
\newcommand {\sing} {\text{sing }}
\newcommand {\ambspace} {\R^N}  
\renewcommand {\max}{{\text{max}}}

\begin{document}

\title{Local Geometry of Singular Real Analytic Surfaces}
\author{Daniel Grieser}
\thanks{The author gratefully acknowledges support
 by the Deutsche Forschungsgemeinschaft (Gerhard-Hess-Programm)
and the Erwin-Schr\"odinger Institute}
\date{January 1999}
\address{Author's address:\\
Institut f\"ur Mathematik\\
Humboldt Universit\"at zu Berlin\\
Unter den Linden 6\\
10099 Berlin\\
Germany}
\begin{abstract}
Let $V\subset\ambspace$ be a compact real analytic surface with isolated singularities,
and assume its smooth part $V_0$ is equipped with a Riemannian metric
that is induced from some analytic Riemannian metric on $\ambspace$. We prove:
\begin{enumerate}
\item Each point of $V$ has a neighborhood which is quasi-isometric (naturally
and 'almost isometrically')
to a union of metric cones and horns, glued at their tips.
\item A full asymptotic expansion, for any $p\in V$,
of the length of $V\cap\{q:\dist(q,p)=r\}$ as $r\to0$.
\item A Gauss-Bonnet Theorem, saying that horns do not contribute an extra
term, while cones contribute the leading coefficient in the length
expansion of 2.
\item The $L^2$ Stokes Theorem, 
self-adjointness and discreteness of the Laplace-Beltrami operator on $V_0$, and
a Gauss-Bonnet Theorem for the $L^2$ Euler characteristic.
\end{enumerate}
As a central tool we use resolution of singularities.
\end{abstract}

\maketitle

\section{Introduction}
In this paper we investigate in detail some aspects of the
geometry and analysis of a compact real analytic surface with isolated
singularities. By this we mean a compact subset $V$ of some 
$\ambspace$ which is locally given as zero set of a finite number 
of real analytic functions, and such that $V$ is a smooth 2-dimensional
manifold outside a finite set $\sing V$. Let $V_0 = V\setminus \sing V$
be the smooth part. 
In order to speak of geometry or analysis on $V$, we need to introduce a
metric. We assume that an {\em induced metric} $g$ is given on $V$, i.e.\
a Riemannian metric on $V_0$ which is induced by some real analytic
Riemannian metric
on $\ambspace$. $g$ induces an (intrinsic) distance function $d$ on $V$. 

Our first concern is a description of the local structure of the
metric space $(V,d)$ up to quasi-isometry. Here, two metric spaces
$(X,d)$, $(X',d')$ are called {\em quasi-isometric} if there is a bijective
map $\Phi:X\to X'$ which is a quasi-isometry, i.e.\ such that
there is a constant
$C>0$ with
$$ C^{-1} d(x,y) \leq d(\Phi(x),\Phi(y)) \leq C d(x,y)\quad\text{ for all }
   x,y\in X.$$
The model metric spaces that will occur are cones and horns. These are 
defined as follows: Let $\gamma\geq 1$ and $\eps>0$. The space
$$ C_\gamma = ([0,\eps)_r\times S^1_\theta)/(\{0\}\times S^1)$$
with distance function $d_\gamma$ induced by the Riemannian metric
$$dr^2 + r^{2\gamma} d\theta^2$$
for $r>0$ and $d_\gamma(0,(r,\theta))=r$, where $0:= [\{0\}\times S^1]$ 
is the 'tip',
is called $\gamma$-horn and $\gamma$ its exponent. A 1-horn is also called a cone.

\begin{theorem}{tqi}
Let $(V,d)$ be as above. Then each $p\in V$ has a neighborhood (in $V$)
which is quasi-isometric to a union of finitely
many cones and horns with rational exponents, glued at their tips.
\end{theorem}

The distance function $d$
on the space glued from $C_{\gamma_1},\ldots,C_{\gamma_n}$
is given by the $d_{\gamma_i}$ on each $C_{\gamma_i}$ and 
$d(p,q)=d_{\gamma_i}(p,0)+d_{\gamma_j}(0,q)$ if $p\in C_{\gamma_i}, q\in C_{\gamma_j}, i\not=j$.

\example If $a\leq b$ are positive integers then the set
$$\{(x,y,z):\, (x^2+y^2)^a = z^{2b} \} \subset \R^3$$
with metric induced from the Euclidean metric on $\R^3$ is quasi-isometric
(near zero)
to two $b/a$-horns glued at their tip.
From this example one sees that any combination of cones and horns can
occur if one allows $V$ to be {\em semi-}analytic (i.e. allowing inequalities
in a local description of $V$) or even semi-algebraic,
so Theorem \ref{tqi}\ gives a complete classification in this case.

Note that a small ball around a smooth point is quasi-isometric to a 
Euclidean ball, which is a cone in polar coordinates.

Theorem \ref{tqi}\ can be sharpened considerably: First, the quasi-isometry
can be chosen to be 'natural', and second, it gets closer and closer to
an isometry as one approaches $p$. Before we state this precisely, we introduce
the infinitesimal notion of quasi-isometry:
Two Riemannian metrics $g,g'$ on a manifold $M$ are called quasi-isometric
if there is $C>0$ such that for each $x\in M$ one has the inequality
of quadratic forms
$$C^{-1} g_x \leq g_x' \leq Cg_x.$$
If this holds, then we write
$$g\qi_C g'$$
or simply $g\qi g'$. $C$ is called the quasi-isometry constant.
It is easy to see that, for two connected Riemannian manifolds $(M,g)$, $(N,h)$, a
diffeomorphism $\Phi:M\to N$ is a quasi-isometry for the induced distance
functions if and only if $\Phi^*h \qi g$.

Let dist denote the (extrinsic) distance function on $\R^N$ and 
$ B_\eps(p) = \{x\in V:\text{ dist}(x,p)<\eps\}.$

\begin{theorem}{tqiexact}
Let $(V,g)$ be as above and let $p\in V$. If $\eps_0$ is sufficiently small
then the following holds for each connected component $K$ of 
$B_{\eps_0}(p) \setminus p$:

There is a  parametrization 
$$\Phi: (0,\eps_0) \times S^1 \to K$$ which 'preserves' $r$, i.e.
\begin{equation} \mylabel{eqrpres}
\dist(\Phi(r,\theta),p)=r \quad\text{ for all }r<\eps_0
\end{equation}
and is normalized
arc length along each $K_r:=\{q\in K: \dist (p,q)=r\}$ such that
\begin{equation} \mylabel{eqquasiiso}
 \Phi^*g \qi_{1+\delta} dr^2 + \left(\frac{l(r)}{2\pi}\right)^2 d\theta^2.
\end{equation} 
Here, $l(r)$ is the length of $K_r$.
In (\ref{eqquasiiso}), one can choose
$$\delta = C\eps^\alpha$$
when restricting $\Phi$ to $(0,\eps)\times S^1, \eps\leq\eps_0$, with 
$C$ and $\alpha>0$ only depending on $V$ and $g$.
\end{theorem}
Explicitly, the arc length condition on $\Phi$ means
\begin{equation}
  |\partial\Phi/\partial\theta|_g = l(r)/2\pi. \mylabel{eqPhinorm}
\end{equation}
Clearly, a parametrization $\Phi$ satisfying (\ref{eqrpres}) and (\ref{eqPhinorm})
(if it exists)
 is uniquely determined by one curve $r\mapsto \Phi(r,\theta_0)$
and an orientation on $K$.

In order to obtain Theorem \ref{tqi}\ from Theorem \ref{tqiexact}\ one
needs to know the leading asymptotic behavior of $l(r)$ as $r\to 0$.
The following theorem gives the complete asymptotics.

\begin{theorem}{tasymp}
The function $l(r)$ defined in Theorem \ref{tqiexact}\
has an asymptotic expansion, as $r\to 0$,
\begin{equation} \mylabel{eqas}
l(r) \sim \sum_{i,j} C_{i,j} r^i(\log r)^j
\end{equation}
with $i$ ranging over $n/m\geq 1$ for some fixed integer $m$ and arbitrary
integers $n\geq m$, and $j$ ranging over $\{0,1\}$. Also, the leading term
contains no logarithm, i.e.\ if 
$$\gamma = \min \{i : C_{i,0}\not=0 \text{ or } C_{i,1} \not=0\}$$
then $C_{\gamma,1}=0$.\\
The asymptotic expansion (\ref{eqas}) may be differentiated (arbitrarily
often) term by
term.
\end{theorem}

We will use these precise descriptions of the local geometry to
draw some conclusions of a global nature. First, we generalize the
classical Gauss-Bonnet theorem.

\begin{theorem}{tgaussbonnet} {\bf (Gauss-Bonnet Theorem)}
Let $V$ be a compact real analytic surface with isolated singularities
$p_1,\ldots,p_R$ with an induced Riemannian metric on 
$V_0=V\setminus\{p_1,\ldots,p_R\}$  as above.
Let $K$ denote
the Gauss curvature.
According to Theorem \ref{tqi}, 
sufficiently small pointed neighborhoods of the $p_i$, $i=1,\ldots,R$,
have connected components that are quasi-isometric
to  cones or  horns. Let $l_1,\ldots,l_S$
be the leading coefficients in the length expansions (Theorem \ref{tasymp})
for all the cone-like components.

Then $K$ is integrable over $V_0$ and
$$ \chi(V) = R +  \frac1{2\pi} \int_{V_0} K  - \frac1{2\pi}\sum_{i=1}^S l_i,
$$
where $\chi(V)$ is the Euler characteristic.
\end{theorem}

We now adress some questions of analysis. The general problem here is the 
following: Given a 
Riemannian manifold $V_0$ which is the smooth part of a compact
'singular space' $V$,
e.g.\ an analytic variety, how do the  properties of the 
natural differential operators defined on $V_0$ differ from those in the smooth
case (i.e.\ $V=V_0$)? Here the 'properties' can be of a very basic nature,
like (i) below, or rather involved, like (iii).
The Riemannian metric $g$ induces a scalar product on differential forms.
By $L^2(\bigwedge V_0)$ we denote the completion of the space of smooth forms
on $V_0$ under this norm. The exterior derivative $d$ and its transpose
$d^t$ act on elements of $L^2(\bigwedge V_0)$ in the sense of distributions.

\begin{theorem}{tspec}
Let $(V,g)$ be as in Theorem \ref{tgaussbonnet}.
\begin{enumerate}
\item[(i)] The $L^2$ Stokes theorem holds, i.e.\ if 
$\omega,\eta,d\omega,d^t\eta\in L^2(\bigwedge V_0)$ then
\begin{equation}\mylabel{eqstokes}
(d\omega,\eta) = (\omega,d^t\eta).
\end{equation}
\item[(ii)] The Gauss-Bonnet operator $\DGB=d+d^t$ and the Laplace-Beltrami
operator $\Delta=\DGB^2$ are self-adjoint as unbounded operators on
$L^2(\bigwedge V_0)$ with
domains $\dom (\DGB) = H^1(\bigwedge V_0)$, $\dom(\Delta) = H^2(\bigwedge V_0)$
(Sobolev spaces). Their spectra are discrete.
\item[(iii)] The $L^2$-Euler characteristic of $V$, i.e.\ the index
of the operator $\DGB$ acting from even forms in $\dom(\DGB)$ to odd forms,
equals
$$\chi_{(2)} (V) = N + \frac1{2\pi} \int_{V_0} K  - \frac1{2\pi}\sum_{i=1}^S l_i,
$$
where $N$ is the total number of cones and horns that occur near singularities
of $V$.
\end{enumerate}
\end{theorem}
(i) is equivalent to a property called 'uniqueness of ideal boundary
conditions', see \cite{GriLes:LSTHTSAV}. It says that no boundary terms occur
in the integration by parts that is implicit in (\ref{eqstokes}).
It is conjectured to be true for all complex projective varieties with
induced metrics, see \cite{CheGorMac:LCIHSAV}\ for related material.
However, it is not true in general for compact real algebraic
varieties, even if the singularities are isolated, see \cite{GriLes:LSTHTSAV}. 
Note that (ii) does not say
that $\Delta$ is essentially self-adjoint on $C_0^\infty(\bigwedge V_0)$. 
In fact, it need not be, as one can see for certain cones.
\vspace{3mm}

\noindent
{\bf Outline of the proofs:}
The asserted quasi-isometry in Theorem \ref{tqiexact}\ with any constant
(resp.\ with the constant $1+O(\eps^\alpha)$) holds iff both of
the following two statements are true (this is elementary; for a proof
see Lemmas \ref{lqi}\ and \ref{lbasicqi}): 
\begin{enumerate}
\item The angle $\beta$ between the curves $\theta=\const$ and the
curves $r=\const$ on $K$ is bounded away from zero as $r\to 0$
(resp. tends to a right angle, with
an error of order $r^\alpha$).
\item The length $|\Phi_r|$ is bounded (resp. is $1+O(r^\alpha)$).
\end{enumerate}
We denote $\Phi_r=\partial\Phi/\partial r$ etc. For simplicity, we
assume in this outline that $p=0$ and $g$ is induced from the Euclidean metric.
Now $\beta=\angle (\Phi_r,\Phi_\theta)$. 
If we write $\Phi(r,\theta) = r\phi(r,\theta)$ with $|\phi|=1$ then 
$\Phi_r = \phi + r\phi_r$ is a decomposition into a 'radial' component $\phi$
which is orthogonal to $\Phi_\theta$ and a 'spherical' component $r\phi_r$.
Therefore, the quasi-isometry is equivalent to the estimate
$r|\phi_r|<C$ (resp. $r|\phi_r|<Cr^\alpha$).

Since the arc length parametrization is difficult to handle directly, we
deal with arbitrary parametrizations $\Psi=r\psi, |\psi|=1$, first and
show (in Section \ref{secqi}) that the $r|\phi_r|$-estimate follows from the analogous
$r|\psi_r|$-estimate plus the uniform boundedness, as $r\to0$,
of the integrals over $K_r$
of the absolute geodesic curvature of $K_r$ in $K$.
The curvature estimate reduces to a length estimate in the unit tangent
bundle, which can be proved by standard analytic techniques, see Section
\ref{seccurv}.

The problem is, then, to construct $\psi$ satisfying $r|\psi_r|<Cr^\alpha$.
This is easy in the model case of a family of curves given by hyperbolas,
i.e. the family $x^ay^b=r$ in the $(x,y)$-plane,
for positive integers $a$ and $b$. Using
resolution of singularities we can reduce the general case to this model.
This will prove Theorem \ref{tqiexact}, see Section \ref{secpf}.

The asymptotic analysis of $l(r)$ (Section \ref{secl}) is also most easily done on
the blown-up space obtained by resolution of singularities. The curves
$K_r$ are then, locally, simply hyperbolas as above, but their length is measured
with respect to a degenerating semi-Riemannian metric.
In such a general situation, one always has an asymptotic expansion as
stated, though the leading term may be $r^\gamma \log r$.
We prove this using the Mellin transform, which reduces the problem
to analytic continuation in the complex variable $z$ of an integral
of the form
$\int_{\R^2} |q(x)|^{z+A} |p(x)|^{1/2} \chi (x) \,dx$
for real analytic functions $p,q$ and $\chi\in C_0^\infty(\R^2)$, plus
a decay estimate for large $\Im z$. To do this, we adapt an
idea by Bernshtejn and Gel'fand \cite{BG}\ and Atiyah \cite{A}, using
resolution of singularities again. 
A separate consideration is needed to exclude a logarithm in the leading
term, in our situation.

The Gauss-Bonnet Theorem is proved (in Section \ref{sgaussbonnet})
by first excising small neighborhoods of
the singular points and applying the classical Gauss-Bonnet Theorem to
the resulting smooth manifold with boundary. The boundary integrals of 
geodesic curvature are interpreted as arc length
variation. Here one needs the differentiated form of (the first term of)
the asymptotics
in Theorem \ref{tasymp}. The integrability of $K$ follows again by a 
standard argument, bounding $|K|$ by the area element of the unit normal bundle
of $V$.

Theorem \ref{tspec}\ follows easily from quasi-isometry
invariance of some of the assertions and well-known results in the
conic and horn-like cases. This is done in Section \ref{sspec}.
\vspace{3mm}

\noindent
{\bf Further remarks and related work:} 
All results are true for semi-analytic sets $V$ as well, since the methods
extend to this case. 
Clearly, the assumption $V\subset\ambspace$ was only made for
notational convenience. All that matters is that the metric on $V_0$
is locally induced by some local embedding. In any case, at least
any smooth real analytic manifold can be embedded analytically in
$\ambspace$ (see \cite{Hir:DT}), so probably the assumption
$V\subset\ambspace$ is no loss of generality.

In the preprint \cite {Bi}, L. Birbrair states a theorem that generalizes
Theorem \ref{tqi}\ to subanalytic surfaces with not necessarily
isolated singularities, for $g$ the Euclidean metric. 
However, there is no analogue of the sharper
Theorem \ref{tqiexact}, and the parametrizations are not explicit. 
J.-M. Lion and J.-P. Rolin
 \cite{LR}\ prove a theorem on asymptotic expansions in a more general
context which implies that the asymptotic series in Theorem \ref{tasymp}\
converges (i.e. $l(r) = l_1(r) + l_2(r)\log r$ with $l_1,l_2$ analytic
in $r^{1/m}$); however, the statement about the leading term is not clear
from \cite {LR}.
The local geometry of semi-algebraic sets is also studied in \cite{BKS}. 
In particular, the authors prove
a Gauss-Bonnet Theorem for semi-algebraic sets. For more on curvature
of sub-analytic sets see J.\ Fu \cite{Fu:CMSS}.

If $V$ is a {\em complex algebraic curve} then it falls in the class of 
surfaces considered here. In this case, Theorem \ref{tqi}\ is much
easier to prove. Questions of analysis have been treated in great detail
by Nagase \cite{Nag:HTSAC},
 Br\"uning-Lesch \cite{BruLes:SGAC},\cite{BruLes:KHTCCC}\
and Br\"uning-Peyerimhoff-Schr\"oder \cite{BruPeySch:OAC}.

The author wishes to thank M.\ Lesch for many helpful discussions and
for suggesting to work on the
spectral theoretic problems that motivated this study,
and
P.\ Milman for clarifying some questions on the
resolution of singularities. Also, thanks to J.\ Fu and J.\ Tolksdorf
for interesting discussions on the subject.


\section{Local metric structure} \mylabel{secpf}
In this section we carry out the main 
steps in the proof of Theorem \ref{tqiexact}.

We will use geodesic normal coordinates around $p$ in $\ambspace$; these are
given by the map
$$ PC: [0,\infty) \times S^{N-1} \to \ambspace $$
defined by $PC(r,\omega) = \gamma(r)$, where $\gamma$ is the geodesic with
$\gamma(0)=p, \gamma'(0)=\omega$. $S^{N-1}$ denotes the unit sphere 
in $T_p\ambspace$
with respect to $g_{|p}$, where we denote by $g$ the metric on $\ambspace$ also.
For $r$ near zero, the form of the metric in normal coordinates is
\begin{equation} \mylabel{eqpolar}
  PC^*g = dr^2 + r^2 \gSN (r)
\end{equation}
where $\gSN(r)$ is a family of Riemannian metrics on $S^{N-1}$ which depends
smoothly on $r\geq 0$. We will denote this metric on $S^{N-1}$ simply
by $|\cdot|$.

We consider local parametrizations of $K$ that satisfy (\ref{eqrpres}),
 i.e.\ which are of  the form
\begin{eqnarray} 
 \Psi: (0,\eps)\times (-A,A) & \to & K, \notag \\
 (r,s) & \mapsto & PC(r,\psi (r,s)) \mylabel{eqparam}
\end{eqnarray}
for some $A>0$, $\eps\leq \eps_0$ and some
map $ \psi : (0,\eps)\times (-A,A) \to S^{N-1}$.

The following proposition reduces the desired quasi-isometry statement
to an estimate on derivatives of parametrizations $\Psi$, and to
curvature estimates. It will be proved
in the next section.

\begin{prop}{pqi}
Suppose that there is a finite collection of parametrizations $\Psi$ as in
(\ref{eqparam})
whose ranges cover $K$ and which
 satisfy the inequality, for some number $C_0$:
\begin{equation}
r|\psi_r|  <  C_0.   \mylabel{eqcond1}
\end{equation}

Also, suppose that there is a constant $C_1$ such that for the geodesic
curvature $\kappa$ of $K_r$ in $K$ we have
\begin{equation} \mylabel{eqgeod}
\int_{K_r} |\kappa|\,d\sigma  < C_1
\end{equation}
where $d\sigma$ denotes arc length.

Then there is a parametrization 
$$\Phi: (0,\eps) \times S^1 \to K$$
satisfying (\ref{eqPhinorm}) and (\ref{eqquasiiso}), with
 $\delta \leq C_0(6+4k(6+2C_1))$
where $k$ is the number of $\Psi$ needed to cover $K$.
\end{prop}

In Section \ref{seccurv}, we prove that (\ref{eqgeod}) is satisfied in
our situation, with $C_1$ only depending on $V$ and $g$.

The following lemma shows that parametrizations $\Psi$ satisfying
a sharper estimate than (\ref{eqcond1}) exist for certain model families
of curves in $\R^2$.

\begin{lemma}{lmodel}
Let $a,b$ be non-negative integers, not both zero, and let $A>0$.
There is a parametrization $\Psi(r,s) = (x(r,s),y(r,s))$
of each connected component of the family of curves
$$ x^ay^b = r,\quad r\in (0,1),\quad |x|+|y|<A, $$
which satisfies
$$ |\Psi_r| \leq r^{\frac1{a+b}-1}.$$
Here, $|\cdot|$ denotes the Euclidean  metric on $\R^2$.
\end{lemma} 

\proof
If $a=0$ then we can take $\Psi(r,s)=(s,\pm r^{1/b})$, and the inequality
is clearly satisfied. Similarly for $b=0$. 

Now assume $a>0,b>0$. We look at the curves in the first quadrant
only, the other components are treated analogously.
The desired parametrization $\Psi=(x,y)$ is defined by setting $s=x-y$. 
Namely, considering
$x,y$ as functions of $r,s$ through the equations $s=x-y, r=x^ay^b$, with
$x,y>0$, one obtains by differentiating these equations with respect to $r$:
$$ x_r = y_r = r^{-1}\left(\frac{a}{x} + \frac{b}{y}\right)^{-1} .$$
By the harmonic-geometric mean inequality, this is bounded above by
$$r^{-1}(a+b)^{-1} (x^ay^b)^{1/(a+b)} = (a+b)^{-1} r^{\frac1{a+b}-1}.$$
The lemma follows.
\qed

In order to reduce the general case to this model, we now
'blow up' $\ambspace$ to resolve the singularity $p$ of $V$. That is, we
find an analytic manifold $X'$ and a proper analytic
map $\pi:X'\to B$, where $B$ is a small ball around $p$ containing
no other singularities of $V$, so that
\begin{enumerate}
\item[(a)] $\pi$ is a diffeomorphism $X'\setminus E \to B\setminus p$
       where $E=\pi^{-1}(p)$.
\item[(b)] The 'strict transform' $V' := \overline{\pi^{-1}(V\setminus p)}$
  is a smooth surface.
\end{enumerate}
Such an (embedded) resolution always exists, see \cite{H}\ or \cite{BM1}.
For our application, we need to know some additional properties of this
resolution:

\begin{lemma}{lres}
The resolution $\pi:X'\to B$ can be chosen so that, in
addition to (a) and (b) above, the following two properties are satisfied:
\begin{enumerate}
\item[(c)] 
  Near any
  point $q\in E\cap V'$ there are coordinates
  $y_1,\ldots,y_N$ centered at $q$ such that, near $q$,
  $$ V' = \{y_3=\ldots=y_N=0\} \quad \text{ and}$$
  $$ r' := \pi^*r = |\prod_{i=1}^N y_i^{a_i}|$$
  for some nonnegative integers $a_i$. Here $r$ denotes the distance
  function to $p$ on $B$.
\item[(d)] Let $h= r^{-2}(g-dr^2)$ be the semi-Riemannian metric on $B\setminus p$
  obtained from the spherical part of $g$. Then there is a smooth semi-Riemannian
  metric $h'$ on $X'$ such that, on $X'\setminus E$, one has $h'=\pi^*h$.
\end{enumerate}
\end{lemma}
\remark Note that $h$ does not extend continuously to $p\in B$. Thus
the point of (d) is that after pulling back to $X'$, a smooth extension
from $X'\setminus E$ to all of $X'$ exists.

Before we prove this lemma, we finish the proof of Theorem \ref{tqiexact}.
By Proposition \ref{pqi}\ and Lemma \ref{lcurv}, 
we are done if we construct parametrizations $\Psi$ covering $V$
which satisfy the estimate (\ref{eqcond1}) with $C_0=C\eps^\alpha$. 
To do this, we first find parametrizations $\Psi'$ of $V'$ and then push them
down to $V$:
For each point $q\in V'_0:=V'\cap E$ let $U_q$ be a neighborhood
of $q$ such that there are coordinates $y_1,\ldots,y_N$ on $U_q$ as in
(c) above. Below we construct, using Lemma \ref{lmodel}, 
parametrizations $\Psi_q'$ of the families
of curves which make up the connected components of  $V'\cap U_q\cap \{r'>0\}$,
such that $\Psi_q'$ satisfies $r|\Psi_{qr}'|_{h'} < C_0=C\eps^\alpha$.
Fix $q$ and set $\Psi=\pi\circ\Psi_q'$. By construction, $\Psi$ 
is of the
form (\ref{eqparam}). We claim that $\psi$ satisfies (\ref{eqcond1}).
Since $(PC^{-1}\Psi)_r = (1,\psi_r)$, we have
$|\Psi_r|_h = |(PC^{-1}\Psi)_r|_\gSN = |\psi_r|_\gSN$. Therefore, we get
$r|\psi_r|_\gSN = r|\Psi_r|_h = r|\Psi'_{qr}|_{h'} \leq C_0$ 
which was to be shown.
Finally, since $r'$ is proper, $V'_0$ is compact,
so finitely many of the $U_q$ cover $V'_0$, and since $r'$ must assume
a positive minimum on the complement of their union $U$, there is an $\eps_0>0$
such that $\{0<r'<\eps_0\}$ is contained in $U$. Thus, the assumptions of
Proposition \ref{pqi}\ are satisfied, so Theorem \ref{tqiexact}\ is proved.

It remains to construct $\Psi_q$ for each connected component of 
$V'\cap U_q\cap \{r'>0\}$. If this set is non-empty then  
$a_3=\ldots =a_N=0$. Then $V'$ is parametrized by $y_1$ and $y_2$ near $q$,
and $r=y_1^{a_1}y_2^{a_2}$. We now use the parametrizations 
given by Lemma \ref{lmodel}.
The estimate (\ref{eqcond1}) is then satisfied
with respect to the Euclidean metric in the $y_1,y_2$-plane and
therefore also with respect to the smooth semi-Riemannian metric $h'$.
Since the exponent of $r$ in Lemma \ref{lmodel}\ is $\frac1{a+b}-1$,
the number $C_0$ can be chosen as $C\eps^\alpha$ whenever $r<\eps$
where $C$ is a
constant only depending on $V$ (stemming from the length distortion
between the Euclidean metric and the metric $h'$ on $U_q$, i.e.
$h'\leq Cg_{\text{eucl}}$ on $U_q$) and 
$$\alpha = \min \frac1{a+b},$$
the minimum being taken over the finitely many coordinate neighborhoods $U_q$.

\qed

\noindent
{\bf Proof of Lemma  \ref{lres}:}

First, we need a few definitions: Let $X$ be a real analytic manifold.
A {\em normal crossings divisor (n.c.d.)} on $X$ is a subset $E\subset X$
which, near any $q\in X$, is a union of coordinate hyperplanes,
in a suitable  local coordinate system near $q$. A set $D$ has normal
crossings with $E$ if the coordinates can be chosen so that, in addition,
$D$ is a coordinate subspace.
Finally, if $f$ is a smooth function on $X$, we say that $f$ has {\em
product form with respect to $E$ (and $D$)} if $\{f=0\}\subset E$ 
and  in one (and therefore any)
such local coordinate
system $y_1,\ldots,y_N$ near $q$ with $y(q)=0$, one has the representation
\begin{equation} \mylabel{eqprod}
 f(y) = h(y)\prod_{i=1}^N y_i^{a_i} 
\end{equation}
with $h$ smooth and $h(0)\not=0$, and non-negative integers $a_i$.

According to \cite {BM1}, a resolution can be obtained in the following
way: 
Define a sequence
$$ \begin{array}{lllllll}
 X_M & \stackrel{\pi_{M-1}}{\to } & X_{M-1} & \to \ldots \to & X_1 &
\stackrel{\pi_0}{\to} & X_0=B \\
E_M && E_{M-1} && E_1 && E_0=\emptyset
\end{array}
$$
of spaces $X_j$, n.c.d.s $E_j\subset X_j$, and proper maps
$\pi_j:X_{j+1}\to X_j$ as follows:
Let $X_0=B, E_0=\emptyset$. If $X_j,E_j$ are defined, choose a suitable
submanifold $D_j\subset X_j$ which is normal crossings with $E_j$,
and let $X_{j+1}$ be the elementary blow-up of  $X_j$ along $D_j$;
that is, essentially $X_{j+1}$ is obtained from $X_j$  by replacing
$D_j$ by its projectivized normal bundle, and $\pi_j$ is the obvious projection
map; also, set
$E_{j+1}=\pi_{j}^{-1}(D_j\cup E_j)$, which is normal crossings again. 
It is one of the main theorems
in \cite{H}\ and \cite{BM1}\ that the $D_j$ can always be found such that
after a finite number of such blowups, $X':=X_M$ is a resolution in the
sense of (a) and (b) above, and such that $V'$ and $E=E_M$
are normal crossings. Here, $\pi=\pi_{(M)}$
where  $\pi_{(j)} = \pi_0\circ\ldots\circ \pi_{j-1}$.
Also, one may choose $D_j$ to lie in the singular set of 
$\overline{\pij^{-1}(V\setminus p)}$ (resp.\ of $V$ for $j=0$).

Assume first that $V$ is singular at $p$.
We claim that (c) and (d) are satisfied for this resolution.
The essential fact is that $V\cap B$ is only singular at $p$,
so that $\pi_0:X_1\to B$ is just the blow-up of $p$. If we identify
$B$ with a subset of $\R^N$ via $x\equiv r\omega$ for $x=PC(r,\omega)$,
 then this can
be described as $X_1= (-\eps_0,\eps_0)\times S^{N-1}/\!\sim$ 
where $(t,\omega)\sim (-t,-\omega)$,
and $\pi([t,\omega]) = t\omega$. Then $E_1=\pi_0^{-1}(0) = \{t=0\}$.
To prove (c), we show inductively that $R_j:=\pij^*r^2$ has product
form with respect to $E_j$, for $j\geq 1$.
This is clearly true for $j=1$ since $\pi_0^*r^2 = t^2$ in the notation above.
Then it follows directly from the definition of elementary blow-up that,
if $R_j$ has product form with respect to $E_j$ and $D_j$ has normal crossings
with $E_j$, then $R_{j+1}=\pi_j^*R_j$ 
has product form with respect to $E_{j+1}$. (Use projective coordinates near 
points in $\pi_j^{-1}(D_j)$.) Therefore, $R_M = \pi^*r^2$ has the form
(\ref{eqprod}). Since $R_M\geq0$,  $h$ is positive and all $a_i$ are even.
Since $E=\{R_M=0\}$, one of the $a_i$ is non-zero for $q\in E$, so one
may modify the coordinates such that $h\equiv 1$. This proves (c).

To prove (d), we note that the pull-back metric on $(-\eps_0,\eps_0)
\times S^{N-1}$ is
$dt^2 + t^2 \gSN(t)$ with $\gSN(t)$ as in (\ref{eqpolar}) for $t\geq0$
and $\gSN(-t)=A^*\gSN(t)$ where $A$ is the antipodal map. This shows
that $h_1 :=\gSN(t)$ is well-defined and smooth on $X_1$,
and satisfies $h_1 = \pi_0^*h$ on $X_1\setminus E_1$.
Therefore, $h'=(\pi_1\circ\ldots\circ \pi_{r-1})^*h_1$ is also smooth.

Finally, if $V$ is smooth at $p$ then just blow up $p$ once, then (c)
and (d) are satisfied by the same argument.
\qed

\section {Basics on quasi-isometry} \mylabel{secqi}

\begin{lemma} {lqi}
Let $\Psi$ be a local parametrization as in (\ref{eqparam}). Then 
\begin{equation}\mylabel{eqqi}
\Psi^*g \qi dr^2 + H(r,s)^2 ds^2
\end{equation}
for some function $H$ if and only if 
$$ r|\psi_r(r,s)| <C $$
for some constant $C$. In this case, one can take $H=r|\psi_s|$.
Then the  quasi-isometry constant  depends only on $C$, and
for $C<1/2$ can be taken to be $1+2C$. 
\end{lemma}

\proof: 
From (\ref{eqpolar}) and $\Psi = PC\circ(r,\psi)$ we have
\begin{eqnarray*}
\Psi^* g &=& dr^2 + r^2 |d\psi|^2 \\
&=& dr^2 + r^2|\psi_r dr + \psi_s ds|^2 \\
&=& (1+|A|^2)dr^2 + 2A\cdot B dr\,ds + |B|^2 ds^2
\end{eqnarray*}
where $A=r\psi_r, B=r\psi_s$.
We now use the following elementary lemma.

\begin{lemma}{lbasicqi}
Let $g=a\,dx^2 + 2b\,dxdy + c\,dy^2$ be a positive definite quadratic
form, and let $\gtilde=a\,dx^2+c\,dy^2$ be its 'diagonal part'.
Set  $$T = \frac{|b|}{\sqrt{ac}}.$$
Then one has
$$ (1-T) \gtilde \leq g \leq (1+T) \gtilde,$$
and these inequalities are sharp.

\end{lemma}
Remark: Note that $T$ is the cosine of the angle between the $x$-axis
and the $y$-axis with respect to the Riemannian metric $g$ on $\R^2$.
Therefore, if $a,b,c$ are allowed to depend on $x$ and $y$, one obtains
that $g\qi \gtilde$ if and only if this angle is bounded away from zero,
and the quasi-isometry constant is close to one if this angle is close
to $\pi/2$.

\proof
The left inequality is equivalent to 
$Ta\,dx^2 + 2b\,dxdy + Tc\,dy^2 \geq 0$, and this
is equivalent to $Ta\geq0$ and $Ta\cdot Tc -b^2 \geq 0$. Similarly for
the right inequality. This implies the claim.
Note that always $T<1$ since $g>0$ implies $ac-b^2>0$.
\qed

Here we have, for fixed $(r,s)$, $a=1+|A|^2, b = A\cdot B$ and $c=|B|^2$,
so $$T=\frac{|A\cdot B|}{|B|\sqrt{1+|A|^2}} \leq \frac{|A|}{\sqrt{1+|A|^2}}.$$
Also, one clearly has
 $dr^2 + |B|^2ds^2 \qi_{1+|A|^2} (1+|A|^2)dr^2 + |B|^2\,ds^2$.
Therefore, from Lemma \ref{lbasicqi}\ we get a quasi-isometry factor, for
fixed $(r,s)$, of at most $(1+|A|^2)(1-\frac{|A|}{\sqrt{1+|A|^2}})^{-1}$.
An easy calculation shows that this is bounded for $A$ bounded, 
and for $|A|<1/2$ is bounded above by $1+2A$. Finally, looking at the length
of $\partial/\partial r$, one sees that $|A|<C$ is also necessary for
the asserted quasi-isometry.

\qed

The following lemma relates arbitrary parametrizations to (unrenormalized)
arc length parametrizations.
\begin{lemma}{lqi2}
Let $\Omega (r,t) = PC(r,\omega(r,t)), (r,t) \in (r_1,r_2)\times (t_1,t_2)$
 be a local parametrization
of $K$ with $\omega$ arc length on the sphere,
i.e.\ $|\omega_t| \equiv 1$. Assume 
that there is a parametrization 
$\Psi:(0,\eps)\times (-A,A)  \to  K$ 
whose range contains the range
of $\Omega$ and which satisfies the inequality (\ref{eqcond1}), and
that (\ref{eqgeod}) holds.

Then, if $r|\omega_r(r,t)|<C$ holds for one value $t=t_0$, it holds for all $t$
(with $C$ replaced by $C+C_0C_2$ where $C_2=6+2C_1$).
\end{lemma}

\proof
We may assume that $\Psi$ and $\Omega$ trace the lines $r=\text{const}$ in
the same direction.
Define the function $b(r)$ by the equation
$$ \omega(r,t_0) = \psi(r,b(r)).$$
Then, if we set $l(r,s) = \int_{b(r)}^s |\psi_s(r,\sigma)| d\sigma$,
we must have 
$$ \omega(r,t_0+l(r,s)) = \psi(r,s).$$
Differentiating in $r$, we get $\omega_r + l_r \omega_t = \psi_r$, evaluated
at corresponding points. 
From condition (\ref{eqcond1}) and $|\omega_t|=1$ we get
\begin{equation}
 r|\omega_r (r,t_0+l(r,s))| < C \Leftrightarrow r|l_r(r,s)| < C,\mylabel{eqequiv}
\end{equation}
where $C$ must be replaced by $C+C_0$ in either direction of the implication.
Now we have 
$l_r = a(r) + \int_{b(r)}^s \frac{\partial}{\partial r}|\psi_s|\, ds$
where $a(r) = -b'(r)|\psi_s(r,b(r))|$. The variation of length formula
says that 
$$ \int_{b(r)}^s \frac\partial{\partial r} |\psi_s| ds = 
\left.(\frac{\psi_s}{|\psi_s|},\psi_r)\right|^s_{b(r)} - \int_{b(r)}^s 
    (N,\psi_r)\kappa\, d\sigma $$
where $\kappa$ is the geodesic curvature, $N$ is the unit normal, and
$\sigma$ denotes arc length. From (\ref{eqcond1}) and (\ref{eqgeod}) we get
$r|l_r - a| \leq C_0(2+C_1)$.
Therefore, the right side of (\ref{eqequiv})
 is equivalent
to $r|a(r)|<C$, where  the constant gets worse by $C_0(2+C_1)$. 
Now since this last condition is independent of $s$, and
since the left side of
(\ref{eqequiv}) is true for $s=b(r)$ by assumption, it must be true
for all $s$. Altogether we lose at most $C_0(6+2C_1)$ in the constant.

\qed

\noindent
{\bf Proof of Proposition \ref{pqi}}

First we consider a parametrization 
\begin{eqnarray*}
 \Omega : (0,\eps) \times [0,\infty) & \to & K\\
         (r,t) & \mapsto & PC(r,\omega(r,t))
\end{eqnarray*}
for which $\omega$ is arc length on the sphere, i.e.\ $|\omega_t| \equiv 1$.
Such a parametrization is determined when  the 'initial' 
curve $r\mapsto\Omega(r,0)$
and an orientation for the family of curves $K_r$ are prescribed. 
We assume such
an orientation given and prescribe the initial curve by choosing any one local
parametrization $\psi$ and setting
\begin{equation} \mylabel{eqinit}
 \omega(r,0) = \psi(r,0).
\end{equation}
Note that $|\Omega_t|=r$ and $\text{length}(K_r)=l(r)$ imply that
\begin{equation}\mylabel{eqperiod}
\omega(r,\lt(r)) = \omega(r,0)
\end{equation}
for $\lt(r) = l(r)/r$.
We now prove that the parametrization $\omega$ thus defined satisfies
the inequality
\begin{equation}\mylabel{eqomegaest}
r|\omega_r(r,t)| < C
\end{equation}
for $t\leq \lt(r)$, with $C=(1+kC_2)C_0$.
Fix $r_0\in (0,\eps)$. We prove (\ref{eqomegaest}) for $r=r_0$, but with
$C$ independent of $r_0$:
Since each of the finitely many given parametrizations $\Psi$ covers an
interval on the curve $\gamma=K_{r_0}$, one can choose points 
$p_0=\Omega(r_0,0),p_1,p_2,\ldots,p_k=p_0$ in this order on $\gamma$ such
that each subarc from $p_i$ to $p_{i+1}$ of $\gamma$ is contained
in the range of some parametrization $\Psi_i$ satisfying the estimates
stated in the proposition. Define $0=t_0<t_1<\ldots<t_k=\lt(r_0)$ by
letting $rt_i$ (for $i<k$) be the distance from $p_0$ to $p_i$ measured along $\gamma$
in the positive direction; then $\Omega(r_0,t_i) = p_i$.  
We prove by induction on $i$ that (\ref{eqomegaest}) is true for 
$t\leq t_i$. For $i=0$, $t_0=0$ and therefore (\ref{eqomegaest}) is true
by (\ref{eqinit}), with $C=C_0$.
Suppose it is true for $t\leq t_i$ with some constant $C=C_{(i)}$.
Since the range of $\Psi_i$ is open, it will contain
the range of $\Omega$ when restricted to a small interval around $r_0$ times
some interval containing $[t_i,t_{i+1}]$. Thus Lemma \ref{lqi2}\ implies
that (\ref{eqomegaest}) holds for $t\in [t_i,t_{i+1}]$, with 
$C=C_{(i+1)} = C_{(i)}+C_2C_0$.
Altogether, we get that (\ref{eqomegaest}) holds for all $t\leq \lt(r)$
with $C=(1+kC_2))C_0$.

Now we renormalize $\omega$ by setting
\begin{equation} \mylabel{eqren}
 \phi (r,\theta) = \omega (r,\theta \,\lt(r)/2\pi).
\end{equation}
We show that $\phi$ also satisfies the inequality
\begin{equation}\mylabel{eqomest}
r|\phi_r(r,\theta)| < C,
\end{equation}
where now $C=(3+2kC_2)C_0$ and $\theta\in [0,2\pi]$.
Differentiating (\ref{eqperiod}) in $r$, we get 
$\omega_r(r,0) = \omega_r(r,\lt(r)) + \lt'(r)\omega_t(r,\lt(r))$.
Now (\ref{eqomegaest}) and $|\omega_t|=1$ give $r|\lt'| < (2+kC_2)C_0$, and 
by (\ref{eqren}) this in turn implies (\ref{eqomest}).

Now set $\Phi=PC\circ (r,\phi)$. By definition, we have $|\Phi_\theta| =
r|\phi_\theta| = r\lt(r)/2\pi = l(r)/2\pi$, and
Lemma \ref{lqi}, applied to $\Phi$,
gives (\ref{eqquasiiso}).
 \qed


\section {Proof of the curvature bound (\ref{eqgeod})} \mylabel{seccurv}

\begin{lemma} {lcurv}
Condition (\ref{eqgeod}) is satisfied, with $C_1$ only depending on $V$ and $g$.
\end{lemma}

\proof
Let $NK_r$ denote the unit normal bundle of $K_r$ in $V$, a one-dimensional
submanifold of  the unit tangent bundle of $V_0$. 
It is an elementary fact from 
differential geometry that
$$ \int_{K_r} |\kappa| \leq \text{length} (NK_r).$$
So it remains to show that these lengths are uniformly bounded as
$r\to 0$.
Set $U=\bigcup_{r\in(0,\eps)} NK_r$. Now $U$ is a semi-analytic subset
of $T^1\ambspace$, the unit tangent bundle of $(\ambspace,g)$. 
This can be seen as follows: Assume $V$ is given as $\{f_1=\ldots=f_L=0\}$
near $p$.
If $V$ is a complete intersection at $p$, 
i.e. $L=N-2$,
then we can write 
$$ U = \{ (x,v)\in T^1\ambspace:
      df_{i|x}(v) = 0, \omega_x(v) = 0, x\in K\setminus p \};$$
here the analytic one-form $\omega$ is defined as
$\omega = *(d(r^2)\wedge df_1\wedge \ldots \wedge df_L)      $
where $*$ denotes the Hodge star operator.
If $V$ is not a complete intersection at $p$ then define, for every subset
$S\subset \{1,\ldots,L\}$ with $N-2$ elements, the one-form
$\omega_S = *(d(r^2)\wedge \bigwedge_{i\in S}df_i)      $
and use the conjunction of $\omega_{S|x}(v)=0$ over all $S$ in the description
of $U$.

We now have the following situation:
$X=T^1\ambspace$ is a Riemannian real analytic manifold, with metric induced
from the metric on $\ambspace$. $U\subset X$ is a two-dimensional
semi-analytic
 subset and $R:=r^2:X\to \R$ is a proper analytic 
function which is non-constant on $U$. We claim that this implies that the level sets
of $R$ have bounded lengths (i.e.\ one-dimensional Hausdorff measure),
for bounded $R$. This can be proved as follows (see \cite{BM2}\ for
definitions and theorems used):
$\Ubar$ is semianalytic and two-dimensional. Let $\pi:Y\to X$ be
a uniformization of $\Ubar$; this means that $Y$ is smooth and two-dimensional,
$\pi$ is proper, and $\pi(Y)=\Ubar$. (A uniformization can be obtained from
a resolution, but also in a simpler way, see \cite{BM2}). Let $h$ be the pull-back semi-Riemannian
metric on $Y$. Now $R'=\pi^*R:Y\to\R$ is proper. Clearly, the length
of any level line $\{R=t\}$ in $X$ is bounded by the length (with respect to
$h$) of the level line $\{R'=t\}$ in Y. Therefore, we are now in the smooth
two-dimensional case, which is easy. For example, one can blow-up
(fairly explicitly) $Y$ to put $R'$ into local product form, then the
level sets are families of hyperbolas locally, and since the pull-back semi-metric is
smooth, they have uniformly bounded length (compare Section \ref{secl}).
\qed


\section {Asymptotics of $l(r)$} \mylabel{secl}
We investigate $l(r)$ by evaluating it on the blown-up space $X'$, described
in Section \ref{secpf}. The 
pull-back metric then degenerates on the preimage of $\{r=0\}$.

We first analyze the model case, where the curves are given in $\R^2$
by $x^ay^b=r$, but the metric may degenerate on the coordinate axes:
\begin{prop}{plasymp}
Let $a,b$ be nonnegative integers, $a+b>0$. Let $h$ be a real analytic
semi-Riemannian
metric defined near zero in $\R^2$ which is Riemannian for $x^ay^b\not=0$, 
and $\chi\in C_0^\infty (\R^2)$ a cutoff function which
is non-negative and equals one near the origin.
Let $l_\chi(r)=l_{\chi,h}(r)$ 
be the length of the curve $x^ay^b=r$, weighted by $\chi$,
i.e.
$$ l_\chi (r) = \int |\psi_s(r,s)|_{h_{\psi(r,s)}} \chi(\psi(r,s)) \, ds $$
for some parametrization $\psi$ of this family of curves.

Then $l_\chi(r)$ has an asymptotic expansion, as $r\to0$
$$l_\chi(r) \sim \sum_{i,j} C_{i,j} r^i(\log r)^j$$
where $i$ ranges over $\frac1m \nnull$ for some fixed
$m\in\N$, and $j\in\{0,1\}$.
This expansion may be differentiated (arbitrarily often) term by term.

If, in addition, $h$ is the pull-back of a non-degenerate (i.e.\
Riemannian) metric under
an analytic map which has injective differential for $x^ay^b\not=0$, 
then the leading term is of the form $r^\gamma$.

\end{prop}

\note The following example shows that, without further assumptions,
 $\log r$ may appear in the leading
term:
Let $h=y^2dx^2 + x^2dy^2$ and $a=b=1$. Then, parametrizing by $s=x$,
one has
$$ |\psi_s|_h = |(1,-r/x^2)|_h = ((r/x)^2 + x^2 (-r/x^2))^{1/2} = \sqrt2\, r/x$$
which is easily seen to imply $l_\chi(r) \sim -\sqrt2\, r\log r$.

\proof
Clearly, $l_\chi$ is bounded, smooth for $r>0$, and has compact support.
We use the Mellin transform, defined, for a bounded
compactly supported function $l$, by
$$ (Ml)(z) = \int_0^\infty r^z l(r) \frac{dr}r$$
with $z\in \C$. This integral is defined and holomorphic for $\Re(z)>0$. 
As is well-known (see \cite{Jea:TMDA}, for example), 
an asymptotic expansion for $l$ as
$r\to 0$ (with derivatives) 
is equivalent to the existence of a meromorphic continuation
of $Ml(z)$ to the whole $z$-plane, with a pole of order $\mu$ at $z_0$
corresponding to a term $r^{-z_0} (\log r)^{\mu-1}$ in the asymptotics, and
such that one has decay for large imaginary part, i.e.
\begin{equation} \mylabel{eqdecay}
|Ml(z)| \leq C_N |\Im z|^{-N}
\end{equation}
for any $N$, uniformly for bounded $\Re z$ and away from the poles.
Therefore, we have to show the existence of such a continuation with at
most double poles, at points in $-\N_0/m$. 
The additional statement on the leading term will
be proved directly.

We assume $a>0,b>0$, the other case is only slightly different, but easier.

A simple calculation, for example using the parameter $s=x$, shows that,
as measures,
$$ |\psi_s|_{h_\psi} \, ds \frac{dr}r = 
                |(bx,-ay)|_{h_{(x,y)}} \frac{dx}x\frac{dy}y.$$
Therefore, we get
\begin{equation} \mylabel{eqml}
 Ml_\chi (z) = \int_0^\infty \int_0^\infty x^{az}y^{bz} p(x,y)^{1/2} \chi(x,y)
               \frac{dx}x\frac{dy}y
\end{equation}               
where we have used $r=x^ay^b$, and where the analytic function $p$ is
defined as $p(x,y)=|(bx,-ay)|^2_{h_{(x,y)}}$. If we set
$x=x_1^2,y=x_2^2$ then  (\ref{eqml}) becomes
$$Ml_\chi(z) =  \int_{-\infty}^\infty \int_{-\infty}^\infty
|x_1|^{2az-1} |x_2|^{2bz-1} p(x_1^2,x_2^2)^{1/2} \chi(x_1^2,x_2^2)
                      \,dx_1 dx_2,$$
and the asymptotics follows from the following proposition:

\begin{prop}{pascont}
Let $p,q$ be real analytic functions defined in a neighborhood $U$
of the support of a function $\rho\in C_0^\infty(\R^n)$. Let $A,B\in\R$,
$B\geq 0$.
For $z\in \C$, $\Re z\geq -A$, define
$$ f(z) = \int_{\R^n} |q(x)|^{z+A} |p(x)|^{B}
           \rho(x) \, dx.$$
Then
\begin{enumerate}
\item[(a)] $f$ can be extended meromorphically to all of $\C$, with poles of 
order at most $n$. 
\item[(b)] If, for every $x\in \supp \rho$,
$q(x)\not=0$ implies $p(x)\not=0$ and $dq_{|x}\not=0$ then
$$ |f(z)| \leq C_N |\Im z|^{-N} $$
for any $N$, uniformly in $|\Im z| \geq 1$ and for bounded $\Re z$.
\end{enumerate}
\end{prop}

\proof
Part (a) is a consequence of a more general theorem by Gel'fand and
Bernshtejn, see \cite{BG}. We sketch their proof and show how one obtains
(b).
First, we use the resolution  of
singularities theorem to find an analytic manifold $X$ and an analytic
proper map $\pi:X\to U$ so that (1) $\pi$ is a diffeomorphism outside 
$\{pq=0\}$ and (2) the pull-back $\pi^*(pq)$ 
has local product form; i.e.\ near any point in $X$ 
there are analytic local coordinates $(y_1,\ldots,y_n)$ centered at that
point such
that $\pi^*(pq)$ is a monomial in the $y_i$ times a non-vanishing analytic 
function.
Then $\pi^*p$ and $\pi^*q$ must have local product form, too, and
if we write $\pi^*(dx_1\wedge\ldots\wedge dx_n) = S(y)dy_1\wedge\ldots\wedge
dy_n$, then $S$ must have local product form (this follows from
 Hilbert's Nullstellensatz -- resp.\ R\"uckert's
Nullstellensatz in the analytic case, see \cite{Rui:BTPS} --
since $S$ vanishes only where $\pi^*(pq)$ vanishes, even after local
complexification).
All this implies that we can use a partition
of unity on $X$ to write $f$ as a finite sum of terms of the form
$$ g(z) = \int_{\R^n} |\qtilde(y)|^{z+A} |\ptilde(y)|^B |S(y)| \rhotilde(y)\,dy$$
where $\rhotilde\in C_0^\infty(\R^n)$, $\qtilde(y)=\pi^*q(y)= y^\alpha h_1(y)$,
$\ptilde(y) = \pi^*p(y) = y^\beta h_2(y)$, $S(y) = y^\gamma h_3(y)$ with
multi-indices $\alpha,\beta,\gamma$ and $h_1,h_2,h_3$ non-vanishing on
$\supp \rhotilde$. For each $g$ we have one of the following two cases:
\begin{enumerate}
\item[(I)] $\alpha\not=0$, i.e.\ $\qtilde(0)=0$.
\item[(II)] $\alpha=0$; then $\gamma=0$ by (1) above, and under
the additional condition in part (b) of the proposition also $\beta=0$,
$d\qtilde\not=0$.
\end{enumerate}
We first consider case (II): If $\alpha=0$ then $|\qtilde|^z$ and therefore
$g$ is entire in $z$, which proves (a). For (b), we may assume w.l.o.g.\ that
$\qtilde(y)=e^{y_1}$, then
$$ g(R+iI) = \int_{-\infty}^\infty e^{iI y_1} \left(
  e^{(R+A)y_1} \int_{\R^{n-1}} |\ptilde|^B |S| \rhotilde\, dy' \right)\, dy_1. $$
Since $\ptilde,S\not=0$ on $\supp \rhotilde$, this is the Fourier transform
of a smooth compactly supported function, therefore rapidly decreasing
in $I$.

We now consider case (I): By a smooth change of coordinates, we may arrange
$h_1\equiv 1$. Then  (with $\delta=\alpha A + \beta B + \gamma$)
$$ g(z) = \int_{\R^n} |y_1|^{\alpha_1 z + \delta_1} \cdots
      |y_n|^{\alpha_n z + \delta_n} \chi(y) \, dy, 
      \quad \chi\in C_0^\infty(\R^n).$$
The following stronger result now implies the proposition, with $m$
equal to the least common multiple of all the {\em positive} $a_i$,
over all $i$ and $g$:
If $\chi\in C_0^\infty(\R^n)$ then the function
$$ g(z_1,\ldots,z_n) = \int_{\R^n} |y_1|^{z_1}\cdots |y_n|^{z_n} \chi(y) \,dy$$
has a meromorphic continuation from $\{\Re z_i\geq 0 \,\forall i\}$ to
$\C^n$, with at most simple poles on the hyperplanes $\{z_i = -s\}, s\in\N$,
and on $\{(z_1,\ldots,z_n):\,\dist(z_i,-\N)>\eps, |\Re z_i|<\eps^{-1} 
\,\forall i\}$
one has for any $N$
$$|g(z_1,\ldots,z_n)| \leq C_{\eps,N} \prod_{i=1}^n |\Im z_i|^{-N}$$
where $C_N$ is bounded in terms of the $C^M$ norm of $\chi$ for some
$M=M(N)$. This claim follows for $n=1$ directly from standard facts about 
distributions, see \cite{GelShi:GF}, and for $n>1$ by induction.
\qed

It remains to prove the additional statement in
Proposition \ref{plasymp} about the leading term.
Assume $h=\pi^*g$ for an analytic map $\pi:(\R^2,0)\to(\R^N,0)$ as in the
proposition and 
a Riemannian metric $g$ on $\R^N$ (everything defined near zero).
Then $|\psi_s|_h = |\partial_s (\pi \circ \psi) |_g$.
Assume $b>0$. We will use the parametrization by $s=x$ now, i.e.
$$ \psi (r,s) = (s,r^{1/b}/s^{a/b}).$$
A short calculation (as above) shows that
\begin{equation} \mylabel{eqdspsi}
 \partial_s (\pi(\psi(r,s))) = \frac1{bx} (bx\pi_x - ay\pi_y) \in \R^N
\end{equation} 
where $(x,y)=\psi(r,s)$. Now the essential fact is the following:
As function of $r$ and $s$, every component of (\ref{eqdspsi})
is a Laurent series in $s^{1/b}$ whose coefficients are
power series in $r^{1/b}$. Since it is a derivative
in $s$, the $s$-residue, i.e.\ the coefficient of $s^{-1}$,
vanishes. This means that no component of the analytic function
$bx\pi_x - ay \pi_y$ contains a monomial of the form $(x^ay^b)^\alpha$.
Also, this function does not vanish identically by the non-degeneracy
hypothesis on $\pi$.
In particular, letting $\gamma\in\Q$ be the largest number such that 
$(x^ay^b)^\gamma$ divides $bx\pi_x-ay\pi_y$, we have
$$bx\pi_x - ay \pi_y = (x^ay^b)^\gamma Q(x,y)$$
with $Q:\R^2\to\R^N$ analytic in certain fractional powers of $x$ and $y$ and
$Q(0,0)=0$, and such that $Q$ does not vanish identially on $\{x=0\}\cup\{y=0\}$.
Since $g$ is nondegenerate, $q=|Q|^2_g$ has the same properties.

To finish the proof, it is clearly enough to 
show that, for some constants $0<c<C<\infty$,
we have $$cr^\gamma < l_\chi(r) < Cr^\gamma,$$
for small $r$. From (\ref{eqdspsi}) and the definition of $q$, we have
$$l_\chi(r) = r^\gamma  I(r)/b \text{ where } I(r) =
                       \int q(x,y)^{1/2} \chi(x,y) \frac{dx}x $$
where always $y=r^{1/b}/x^{a/b}$.
To show $I<C$, we observe that $q(0,0)=0$ implies $q(x,y) \leq x^B+y^B$ for
some $B>0$. Since $\frac{dx}x = \frac{b}a \frac{dy}y$,
we have $I \leq \int_0^1 x^B dx/x + b/a\int_0^1 y^B dy/y < C$.

To show $I>c$, we recall that $q$ does not vanish identically on at least
one coordinate axis, say the $x$-axis. Then, $q^{1/2}>c_0>0$ on some rectangle
$J\times [0,\delta]$ on which $\chi=1$, where $J$ is an open interval
on the $x$-axis and $\delta>0$. Then $I \geq c_0 \int_J dx/x =:c$ whenever $r$ is so 
small that $x^ay^b=r, x\in J$ implies $y<\delta$.
This finishes the proof of the proposition.

\qed

\noindent {\bf Proof of Theorem \ref{tasymp}:}
Apply the proposition to the semi-Riemannian metric $h_1=\pi^* g$
and to $h_2= h'$ from Lemma \ref{lres}\ (d), and use
a partition of unity on $X'$.
Clearly $l_{\chi,h_1} = rl_{\chi,h_2}$. Therefore, the leading term
has exponent $\gamma\geq 1$. Also, since $g$ is non-degenerate,
there is no logarithm.
\qed


\section{The Gauss-Bonnet Theorem} \mylabel{sgaussbonnet}
Here we prove Theorem \ref{tgaussbonnet}.
The integrability of $K$, i.e.\ the finiteness of $\int_{V_0} |K|$,
is proved in a similar way as the bound on geodesic curvature (\ref{eqgeod}),
only simpler. We sketch the argument in three dimensions; for details
see \cite{Eul:IKAURKK}:
From the fact that $K$ is the Jacobian of the differential
of the Gauss map one easily sees that for any open, relatively compact set
$U\subset V_0$ one has $\int_U |K| \leq \text{area}(NU)$, where $NU$ is
the unit normal bundle of $U$, with the induced metric. Since $V_0$ is 
semi-analytic, so is $NV_0$, so it has finite area by compactness
of $V$. Therefore, $\int_{V_0} |K|$ is finite also.

To prove the Gauss-Bonnet formula, we choose $\eps_0>0$ small enough
so that Theorem \ref{tqiexact}\ applies for each $p_i$  
(but here we only need the 
topological statement), and let
$U_\eps = \bigcup_{i=1}^R B_\eps(p_i)$ and $V_\eps = V\setminus U_\eps$.
For $\eps<\eps'<\eps_0$ use the additivity of the Euler characteristic
and $V=V_\eps \cup U_{\eps'}$ to obtain
\begin{equation}\mylabel{eqdecomp}
\chi(V) = \chi(V_\eps) + \chi(U_{\eps'})  - \chi(V_\eps\cap U_{\eps'}).
\end{equation}
Now $U_{\eps'}$ is homotopic to $\{p_1,\ldots,p_R\}$, so 
$$ \chi(U_{\eps'}) = R,$$
and $V_\eps\cap U_{\eps'}$ is homotopic to a disjoint union of circles, so
$$ \chi(V_\eps\cap U_{\eps'}) = 0.$$
Applying the classical Gauss-Bonnet Theorem (see \cite{Spi:CIDG}, Ch.\ 6)
to the manifold with boundary
$V_\eps$, we get
\begin{equation} \mylabel{eqclassgb}
\chi(V_\eps) = \frac1{2\pi}\int_{V_\eps} K + \frac1{2\pi} \int_{\partial V_\eps}\kappa,
\end{equation}
where $\kappa$ is the geodesic curvature of $\partial V_\eps = \partial U_\eps$
with respect to the normal pointing outward from $V_\eps$.
Now the variation of arc length is geodesic curvature
(see \cite{Spi:CIDG}, Ch.\ 9), i.e.\ (with the previous orientation)
$$ \frac{d}{d\eps} l(\partial V_\eps) = 
    - \int_{\partial V_\eps} (\partial/\partial r, N)\kappa$$
where $N$ is the unit normal (pointing away from the $p_i$) and
$\partial/\partial r$ is taken with respect to the coordinates described in
Theorem \ref{tqiexact}. Now this theorem implies that the scalar product
$(\partial/\partial r, N)$ tends to 1 as $r=\eps\to 0$, uniformly in $\theta$.
Since $\int_{V_\eps} |\kappa|$ is uniformly bounded (Lemma \ref{lcurv}), we get
$$ \frac{d}{d\eps} l(\partial V_\eps) = 
    - \int_{\partial V_\eps} \kappa + o(\eps)$$
By Theorem \ref{tasymp}, we have 
$$ \frac{d}{d\eps} l(\partial V_\eps)_{|\eps=0} = \sum_{i=1}^S l_i.$$
Therefore, taking the limit $\eps\to0$ in (\ref{eqdecomp}) and (\ref{eqclassgb})
we get the result.
\qed

\section{Analysis near a singular point} \mylabel{sspec} 
Here we prove Theorem \ref{tspec}.
For a differential operator $D$ on forms, denote by $D_{\min}$ and
$D_{\max}$ the unbounded operators
in $L^2:=L^2(\bigwedge V_0)$ with domains 
$\{\omega\in L^2: \exists \omega_n\in C_0^\infty(V_0),\, \omega_n\to\omega
\text{ and } D\omega_n\to D\omega\text{ in }L^2\}$ and
 $\{\omega\in L^2: D\omega\in L^2\}$, respectively.

Proof of (i): First, we show that validity of the $L^2$ Stokes Theorem
depends only on the quasi-isometry class of the metric. This is not
obvious, since the local expression for $d^t$ contains derivatives of
the metric coefficients. But it follows directly from the fact that
(\ref{eqstokes}) is equivalent to the statement
$$ d_{\min} = d_{\max}$$
and the fact that these two operators (in particular their domains)
are quasi-isometry invariants.
To see this equivalence, note that (\ref{eqstokes}) is equivalent to 
\begin{equation}\label{eqstokes1}
(d_{\max})^* = (d^t)_{\max}
\end{equation}
where the star denotes the functional analytic adjoint. Also, it is
easy to see that $(d_{\min})^* = (d^t)_{\max}$. Since $d_{\min}$ and 
$d_{\max}$ are closed
operators, the claim follows. For details see \cite{GriLes:LSTHTSAV}.

By Theorem \ref{tqi}\ it is therefore enough to prove (\ref{eqstokes})
for metrics on $V_0$ that are exactly horn- or cone-like near the singularities.
This is a special case of a well-known result  by Cheeger \cite{Che:HTRP},
Theorem 2.2.

Proof of (ii): 
By (i) and (\ref{eqstokes1}), the operator $d_\max+(d^t)_\max$ is self-adjoint.
Its domain is $\{\omega\in L^2: d\omega, d^t\omega\in L^2\}=
\dom(D_{\text{GB,max}})
=: H^1(\bigwedge V_0)$, where the first equality follows from the orthogonality of 
the ranges of $d$ and $d^t$.
Since $H^2(\bigwedge V_0) = \{\omega\in L^2: 
d\omega,\Delta\omega\in L^2\}=\dom((D_{\text{GB,max}})^2)$, one also
has self-adjointness of $\Delta$. Finally, discreteness of the spectrum
is a quasi-isometry invariant by \cite{BruLes:HC}, Lemma 2.17, and
was proved for cones and horns by Cheeger \cite{Che:SGSRS}\ and
Lesch-Peyerimhoff \cite{LesPey:IFMMH}.

Proof of (iii):
One can argue as in the proof of Theorem \ref{tgaussbonnet}: $\chi_{(2)}$
is also additive. For manifolds with boundary (and metrics smooth up to
the boundary) or their open interior, one
has $\chi_{(2)} = \chi$, so only the term $\chi_{(2)}(U_{\eps'})$ in
(\ref{eqdecomp}) needs to be reevaluated. Now $\chi_{(2)}$ is a quasi-isometry
invariant, and on cones and horns it equals one by \cite{Che:HTRP}, Lemma 3.4.
Since $\chi_{(2)}$ is calculated on the smooth part $V_0$, one
has $\chi_{(2)}(U_{\eps'}) = \chi_{(2)}(U_{\eps'}\setminus{p_1,\ldots,p_R})
=N$. \qed

\end{document}